\documentclass[12pt,a4paper,reqno]{article}
\usepackage[a4paper,portrait]{geometry}
\usepackage{amsmath,amsfonts,amssymb,amsthm}

\theoremstyle{plain}%
\newtheorem{thm}{Theorem}[section]%
\newtheorem{prop}[thm]{Proposition}%
\newtheorem{lem}[thm]{Lemma}%
\theoremstyle{definition}%
\theoremstyle{remark}%

\newcommand{\EM}{\ensuremath}

\newcommand{\dH}{\EM{\mathbb{H}}}

\newcommand{\dN}{\EM{\mathbb{N}}}

\newcommand{\dR}{\EM{\mathbb{R}}}

\newcommand{\dZ}{\EM{\mathbb{Z}}}

\newcommand{\cC}{\EM{\mathcal{C}}}

\newcommand{\cH}{\EM{\mathcal{H}}}

\newcommand{\cL}{\EM{\mathcal{L}}}

\newcommand{\cO}{\EM{\mathcal{O}}}
\newcommand{\cP}{\EM{\mathcal{P}}}

\newcommand{\cZ}{\EM{\mathcal{Z}}}

\newcommand{\bE}{\EM{\mathbf{E}}}

\newcommand{\bI}{\EM{\mathbf{I}}}

\newcommand{\bL}{\EM{\mathbf{L}}}

\newcommand{\al}{\alpha}
\newcommand{\be}{\beta}

\newcommand{\de}{\delta}
\newcommand{\ga}{\gamma}
\newcommand{\Ga}{\Gamma}

\newcommand{\la}{\lambda}
\newcommand{\La}{\Lambda}

\newcommand{\si}{\sigma}

\newcommand{\veps}{\varepsilon}
\newcommand{\vphi}{\varphi}

\newcommand{\p}[4]{{#3}\!\left#1{#4}\right#2} 

\newcommand{\ABS}[1]{\EM{{\left| #1 \right|}}} 
\newcommand{\BRA}[1]{\EM{{\left\{#1\right\}}}} 
\newcommand{\DP}[1]{\EM{{\left<#1\right>}}} 
\newcommand{\NI}[1]{\EM{{\left\| #1\right\|}_\infty}} 
\newcommand{\OSC}[1]{\EM{\mathrm{osc}}\PAR{#1}} 
\newcommand{\PAR}[1]{\EM{{\left(#1\right)}}} 
\newcommand{\BPAR}[1]{\EM{{\biggl(#1\biggr)}}} 
\newcommand{\pd}{\EM{\partial}} 
\newcommand{\SBRA}[1]{\EM{{\left[#1\right]}}} 

\newcommand{\entf}[1]{\mathbf{Ent}_{#1}}
\newcommand{\ent}[2]{\p(){\entf{#1}}{#2}}

\newcommand{\moyf}[1]{\bE_{#1}}
\newcommand{\moy}[2]{\p(){\moyf{#1}}{#2}}
\newcommand{\bmoy}[2]{\moyf{#1}\!\biggl(#2\biggr)}

\newcommand{\covf}[1]{\mathbf{Cov}_{#1}}
\newcommand{\cov}[3]{\p(){\covf{#1}}{#2,#3}}
\newcommand{\bcov}[3]{\covf{#1}\!\biggl(#2,#3\biggr)}

\newcommand{\varf}[1]{\mathbf{Var}_{#1}}
\newcommand{\var}[2]{\p(){\varf{#1}}{#2}}
\newcommand{\bvar}[2]{\varf{#1}\!\biggl(#2\biggr)}

\newcommand{\Det}[1]{\mathrm{Det}\,}

\newcommand{\Mo}{\mathbf{1}}
\newcommand{\Vo}{\mathrm{1}}

\newcommand{\GI}{\bL}
\newcommand{\GR}{\nabla}

\newcommand{\LA}{\boldsymbol{\Delta}}

\newcommand{\Hess}[1]{{\p(){\mathrm{Hess}}{#1}}}

\newcommand{\SSK}[1]{\substack{#1}}
\renewcommand{\leq}{\leqslant}
\renewcommand{\geq}{\geqslant}
\newcommand{\bs}{\EM{\backslash}}

\newcommand{\mykeywords}{Interacting particle systems, spectral gap,
  Poincar\'e inequality, Log-Sobolev inequality, Conservative spin systems,
  Continuous spin systems, Ginzburg-Landau process on a lattice, Glauber
  Dynamics, Kawasaki Dynamics, Mean-field models, Exchangeable measures}
        
\newcommand{\mysubjclass}{
60K35 Interacting random processes, 
82B44 Disordered systems, 
82B20 Lattice systems,
46-99 Functional analysis,
60J60 Diffusion processes,
26D10 Inequalities involving derivatives, differential and integral operators}

\title{Glauber versus Kawasaki for spectral gap\\ and logarithmic Sobolev
  inequalities\\ of some unbounded conservative spin systems}

\date{University of Toulouse\\ \medskip \small November 5, 2002\\%
      \footnotesize \medskip Compiled \today}

\author{D.~Chafa\"{\i}}

\newcommand{\MSX}{\EM{M-\sum_{i=1}^n x_i}}
\newcommand{\muMi}{\EM{\mu_{(i)}}}
\newcommand{\siMi}{\EM{\si_{(i)}}}

\newcommand{\sik}[1]{\EM{\si^{(#1)}}}  

\newcommand{\LAW}[2]{\EM{\mathrm{Law}_{#1}\PAR{#2}}}

\begin{document}

\maketitle

\begin{abstract}
  Inspired by the recent results of C. Landim, G. Panizo and H.-T. Yau
  \cite{landim-panizo-yau} on spectral gap and logarithmic Sobolev
  inequalities for unbounded conservative spin systems, we study uniform
  bounds in these inequalities for Glauber dynamics of Hamiltonian of the form
  $$
  \sum_{i=1}^n V(x_i)+V\biggl(\MSX\biggr), \quad (x_1,\ldots,x_n)\in\dR^n
  $$
  Specifically, we examine the case $V$ is strictly convex (or
  small perturbation of strictly convex) and, following
  \cite{landim-panizo-yau}, the case $V$ is a bounded perturbation of
  a quadratic potential. By a simple path counting argument for the
  standard random walk, uniform bounds for the Glauber dynamics
  yields, in a transparent way, the classical $L^{-2}$ decay for the
  Kawasaki dynamics on $d$-dimensional cubes of length $L$. The
  arguments of proofs however closely follow and make heavy use of
  the conservative approach and estimates of \cite{landim-panizo-yau},
  relying in particular on the Lu-Yau martingale decomposition and
  clever partitionings of the conditional measure.
\end{abstract}

\section*{Introduction}

Let $Q$ be a probability measure on $\dR^n$. In the sequel, we denote
by $\moy{Q}{f}$ the expectation of $f$ with respect to $Q$,
$\var{Q}{f}:=\moy{Q}{f^2}-\moy{Q}{f}^2$ the variance of $f$ for $Q$,
and $\ent{Q}{f}$ the entropy of a non negative measurable function $f$
with respect to $Q$, defined by
$$
\ent{Q}{f}:=\int\! f\log f dQ - \int\! f dQ \;\log \int\! f dQ.
$$
We say that $Q$ satisfies a Poincar\'e inequality if there exists a positive
constant $\cP$ such that for any smooth function $f:\dR^n\to\dR$,
\begin{equation}\label{in:def-ts}
\var{Q}{f} \leq \cP\, \moy{Q}{\ABS{\GR f}^2},
\end{equation}
where $\ABS{\GR f}^2:=\sum_{i=1}^n \ABS{\pd_i f}^2$. Similarly, we say that
$Q$ satisfies a logarithmic Sobolev inequality if there exists a positive
constant $\cL$ such that for any smooth function $f$,
\begin{equation}\label{in:def-ls}
\ent{Q}{f^2} \leq \cL\, \moy{Q}{\ABS{\GR f}^2}.
\end{equation}
This inequality strengthen the Poincar\'e inequality \eqref{in:def-ts} since for
$\veps$ small enough,
$$
\ent{Q}{(1+\veps f)^2}=2\veps^2\var{Q}{f}+\cO(\veps^3),
$$
which gives $2\cP \leq \cL$. Let $\cH\in\cC^2(\dR^n,\dR)$ such that
$$
Z_{\cH}:=\int_{\dR^n}\!\!e^{-\cH(x)}\,dx < +\infty.
$$
The probability measure $Q$ defined by
$dQ(x)=(Z_{\cH})^{-1}\exp\PAR{-\cH(x)}\,dx$ is the symmetric invariant measure
of the diffusion process $(X_t)_{t\geq 0}$ on $\dR^n$ driven by the S.D.E.
$$
dX_t = \sqrt{2}\,dB_t - \GR \cH(X_t)\,dt,
$$
where $(B_t)_{t\geq 0}$ is a standard Brownian motion on $\dR^n$. In this
context, we say that the probability measure $Q$ is associated with the
``Hamiltonian'' $\cH$. It is well known that $Q$ satisfies the Poincar\'e
inequality \eqref{in:def-ts} with a constant $\cP$ if and only if the
infinitesimal generator $\GI:=\LA-\GR\cH\cdot\GR$ possesses a spectral gap
greater than $\cP^{-1}$. In the other hand, a famous Theorem of Gross states
that $Q$ satisfies the logarithmic Sobolev inequality \eqref{in:def-ls} if and
only if the diffusion semi-group generated by $\GI$ is hyper-contractive. A
celebrated result of Bakry and \'Emery ensures that when there exists a
constant $\rho>0$ such that for any $x\in\dR^n$,
$$
\Hess{\cH}(x)\geq \rho\,\bI_{p}
$$
as quadratic forms on $\dR^n$, i.e. $\cH$ is uniformly strictly convex or
$Q$ is log-concave, then $Q$ satisfies to \eqref{in:def-ts} and
\eqref{in:def-ls} with constants $\cP=\rho^{-1}$ and $\cL=2\,\rho^{-1}$
respectively. Moreover, $Q$ satisfies to \eqref{in:def-ts} with a constant
$\cP$ if and only if
$$
\cP\moy{Q}{(\GI f)^2} \geq \moy{Q}{\ABS{\GR f}^2}
$$
for any smooth function $f$. The reader may find an introduction to
logarithmic Sobolev inequalities and related fields in \cite{our-ls-2000}.
 
We are interested in the present work to particular ``Hamiltonians''
$H\in\cC^2(\dR^{n+1},\dR)$. Let $M\in\dR$ and define $H_M\in\cC^2(\dR^n,\dR)$
by
\begin{equation}\label{eq:defhm}
  H_M(x_1,\ldots,x_n):=H\biggl(x_1,\ldots,x_n,\MSX\biggr).
\end{equation}
Assume that $Z_{H_M}<\infty$. Our aim is to establish Poincar\'e and
logarithmic Sobolev inequalities for probability measures on $\dR^n$ of the
form
\begin{equation}\label{eq:defnum}
d\si_M(x_1,\ldots,x_n):=(Z_{H_M})^{-1}\,\exp\PAR{-H_M(x)}\,dx_1\cdots dx_n,
\end{equation}
with constants $\cP$ and $\cL$ which does not depend on $n$ and $M$.  This
investigation is motivated by the study of certain conditional probability
measures. Namely, if the probability measure $\mu$ on $\dR^{n+1}$ given by
\begin{equation}\label{eq:defmu}
d\mu(x):=(Z_H)^{-1}\,\exp\PAR{-H(x_1,\ldots,x_{n+1})}\,dx_1\cdots dx_{n+1}
\end{equation}
is well-defined, i.e. $Z_H<+\infty$, then for any $M\in \dR$, one can define
the conditional probability measure $\mu_M$ by
\begin{equation}\label{eq:defmum}
\mu_M:=\mu\biggl(\;\cdot\;\Bigr\vert\;\sum_{i=1}^{n+1} x_i=M\biggr),
\end{equation}
and we get, for any $f\in\cC_b(\dR^{n+1},\dR)$,
\begin{equation}\label{eq:mum2num} 
  \moy{\mu_M}{f}=\int_{\dR^n}\! f\biggl(x_1,\ldots,x_n,\MSX\biggr)\,
  d\si_M(x_1,\ldots,x_n). 
\end{equation}
Thus, $\si_M$ can be viewed as the translation of the conditional
probability measure $\mu_M$ under the affine hyper-plane of $\dR^{n+1}$
of equation $x_1+\cdots+x_{n+1}=M$. Alternatively, and following
Caputo in \cite{caputo}, the conditional probability measure $\mu_M$
can be defined from the probability measure $\mu$ given in
\eqref{eq:defmu} by adding an infinite potential outside of the affine
constraint $x_1+\cdots+x_{n+1}=M$. Namely for any bounded continuous
function $f:\dR^{n+1}\to\dR$
$$
\moy{\mu_M}{f}=\lim_{\be\to+\infty} \moy{\mu_{M,\be}}{f}.
$$
where $\mu_{M,\be}$ denotes the probability measure on $\dR^{n+1}$
defined by
$$
d\mu_{M,\be}(x) 
:= \cZ_{\mu_{M,\be}}^{-1}\,e^{-\be\,\PAR{M-x_1-\cdots-x_{n+1}}^2}\,d\mu(x).
$$
A simple change of variable in $\moy{\mu_{M,\be}}{f}$ gives that
$$
\lim_{\be\to+\infty} \moy{\mu_{M,\be}}{f}=
\int_{\dR^n}\! f\biggl(x_1,\ldots,x_n,\MSX\biggr)\,d\si_M(x_1,\ldots,x_n),
$$
This weak limit definition of $\mu_M$ was used by Caputo in
\cite{caputo} in order to study the case of a convex Hamiltonian $H$.
We do not use it in our approach.  Notice that if
$f\in\cC_b(\dR^n,\dR)$, we get from \eqref{eq:mum2num} that
\begin{equation}\label{eq:mumeqnum} 
  \moy{\si_M}{f} 
  =\int_{\dR^{n+1}}\!f(x_1,\ldots,x_n)\, d\mu_M(x_1,\ldots,x_{n+1}).
\end{equation}
Observe that \eqref{eq:mum2num} gives $\moy{\mu_M}{x_1+\cdots+x_{n+1}}=M$.
Thus, when $H$ is a symmetric function, $\si_M$ and $\mu_M$ are exchangeable
measures, i.e. invariant by any permutation of the coordinates. This holds for
example when $H(x)=V(x_1)+\cdots+V(x_{n+1})$.  Moreover, $\MSX$ and $x_j$ have
then the same law under $\si_M$ for any $j$ in $\{1,\ldots,n\}$ and we get
\begin{equation}\label{eq:mum-num-moy}
  \moy{\mu_M}{x_1}=\cdots=\moy{\mu_M}{x_{n+1}}
  =\moy{\si_M}{x_1}=\cdots=\moy{\si_M}{x_n}
  =\frac{M}{n+1}.
\end{equation}
Thus, the mean of $\mu_M$ and $\si_M$ does not depend on $H$ in this case.

Let us see now how to translate \eqref{in:def-ts} and \eqref{in:def-ls} for
$\si_M$ in terms of $\mu_M$. One can observe that for any $i\in\{1,\ldots,n\}$
\begin{multline*}
  \pd_i\BPAR{f\BPAR{x_1,\ldots,x_n,\MSX}} \\
  =(\pd_i f)\BPAR{x_1,\ldots,x_n,\MSX}-(\pd_{n+1} f)\BPAR{x_1,\ldots,x_n,\MSX}.
\end{multline*}
By replacing the coordinate $x_{n+1}$ by any of the $x_1,\ldots,x_n$ in
\eqref{eq:mum2num}, we obtain the following proposition
\begin{prop}\label{pr:num-mum}
  Let $H:\dR^{n+1}\to\dR$ and assume that for any permutation $\pi$ of the
  coordinates, the probability measure $\si_M^\pi$ on $\dR^n$ defined by
  \eqref{eq:defnum} and associated to $H\circ \pi$ satisfies to Poincar\'e
  (resp.  logarithmic Sobolev) inequality with a constant $\cP$ (resp. $\cL$)
  which does not depend on $n$, $M$ and $\pi$.  Then, if $\mu_M$ is the
  associated conditional probability measure defined by \eqref{eq:defmum}, we
  get for any smooth $f:\dR^{n+1}\to\dR$
  \begin{equation}\label{in:ts-ij}
   \var{\mu_M}{f} 
   \leq 
   \frac{\cP}{n+1}
   \;\moy{\mu_M}{\sum_{1\leq i,j\leq n+1}\ABS{\pd_i f-\pd_j f}^2},  
  \end{equation}
  and respectively
  \begin{equation}\label{in:ls-ij}
  \ent{\mu_M}{f^2} 
  \leq 
  \frac{\cL}{n+1}
  \;\moy{\mu_M}{\sum_{1\leq i,j\leq n+1}\ABS{\pd_i f-\pd_j f}^2}.
\end{equation}
\end{prop}
These inequalities leads to constants in $L^2$ for the ``Kawasaki
dynamics'' associated to $\mu_M$. Namely, consider a finite box
$\La:=\{1,\ldots,L\}^d\subset\dZ^d$ on the lattice $\dZ^d$ and $n$
such that $\dR^\La\simeq\dR^{n+1}$ (i.e. $n+1=\ABS{\La}=L^d$). There
exists a constant $C>0$ depending only on $d$ such that for any
$a\in\dR^\La$
\begin{equation}\label{in:compdirf} 
 \frac{1}{\ABS{\La}}\sum_{i,j\in\La}(a_i-a_j)^2
 \leq 
 C\;L^2\;\sum_{\SSK{i,j\in\La\\\ABS{i-j}=1}}(a_i-a_j)^2,
\end{equation}
Therefore, it is straightforward to deduce from \eqref{in:ts-ij} and
\eqref{in:ls-ij} that for a constant $C>0$ which does not depend on
$n$ and $M$, one have
\begin{equation}\label{in:ts-l2}
\var{\mu_M}{f} 
\leq C\,L^2\; \sum_{\SSK{k,l\in\La\\\ABS{k-l}=1}}
\moy{\mu_M}{\ABS{\pd_i f-\pd_j f}^2},
\end{equation}
and
\begin{equation}\label{in:ls-l2}
\ent{\mu_M}{f^2} \leq
C\,L^2\;
\sum_{\SSK{k,l\in\La\\\ABS{k-l}=1}}
\moy{\mu_M}{\ABS{\pd_i f-\pd_j f}^2}.
\end{equation}
Inequality \eqref{in:compdirf} follows from a classical path counting argument
(see for example section 4.2 of \cite{saloffcoste}). However, let us gives
briefly a proof. For any $i,j$ in $\La$, consider the path $\Ga_{ij}$ inside
$\La$ joining $i$ and $j$ obtained by adjusting the $d$ coordinates one after
the other. We have $\ABS{\Ga_{ij}}\leq dL$ and for each $k,l$ in $\La$ with
$\ABS{k-l}=1$, the number of such paths containing the edge $(k,l)$ is bounded
above by $c_d\,L^{d+1}$ where $c_d>0$ is a constant depending only on $d$. Now
by Cauchy-Schwarz's inequality
$$
(a_i-a_j)^2
=
\SBRA{\sum_{(k,l) \in \Ga_{i,j},\,\ABS{k-l}=1}\!\!(a_k-a_l)}^2
\leq
dL \sum_{(k,l) \in \Ga_{i,j},\,\ABS{k-l}=1} (a_k-a_l)^2,
$$
and therefore
$$
\sum_{i,j\in\La} (a_i-a_j)^2
\leq
dL \sum_{\SSK{\ABS{k-l}=1\\k,l\in \La}}
(a_k-a_l)^2\sum_{\SSK{i,j\in\La\\\Ga_{i,j}\ni(k,l)}}\!\!\!1
\leq 
dc_d L^{d+2} \sum_{\SSK{\ABS{k-l}=1\\k,l\in \La}} (a_k-a_l)^2,
$$
which gives the desired result \eqref{in:compdirf}.

A simple example is given by uniformly strictly convex $H$ in $\dR^{n+1}$.
Namely, if there exists a constant $\rho>0$ such that for any $x\in\dR^{n+1}$,
$\Hess{H}(x)\geq \rho \bI_{n+1}$ as quadratic forms on $\dR^{n+1}$, then, an
easy calculus gives for any $x\in\dR^n$ and $h\in\dR^n$
$$
\DP{(\Hess{H_M})(x) h,h}_{\dR^n}
\geq \rho
\sum_{i=1}^n h_i^2 + \rho \PAR{-\sum_{i=1}^n h_i}^2 \geq \rho \sum_{i=1}^n h_i^2.
$$
Thus, $H_M$ is uniformly strictly convex with the same constant
$\rho$, and therefore, by the Bakry-\'Emery criterion, $\si_M$ satisfies
to Poincar\'e and logarithmic Sobolev inequalities with a constant
$\rho^{-1}$ and $2\,\rho^{-1}$ respectively, which does not depend on
$n$ and $M$.  The hypotheses of Proposition \ref{pr:num-mum} are
full-filled since by the same calculus, $(H\circ\pi)_M$ is also
uniformly strictly convex with a constant $\rho$. A more simple
example is given by 
$$
H(x)=V(x_1)+\cdots+V(x_{n+1})
$$
where $V$ is in $\cC^2(\dR,\dR)$ with $V''>\rho>0$.  Let us
consider now another convex Hamiltonian example on $\dR^{n+1}$ defined
by
$$
H(x):=\frac{1}{2(n+1)}\,\sum_{i,j=1}^{n+1} V_{\{i,j\}}(x_i-x_j),
$$
where $V_{\{i,j\}}$ are in $\cC^2(\dR,\dR)$ and even.  This is a so called
mean-field Hamiltonian when all the $V_{\{i,j\}}$ are equal.  We have for any
$i,j$ in $\{1,\ldots,n+1\}$
$$
(n+1)\,\pd^2_{ij} H(x)=
\begin{cases}
  \sum_{\SSK{k=1\\k\neq i}}^{n+1} V''_{\{i,k\}}(x_i-x_k) & \text{ if } i=j\\
  -V''_{\{i,j\}}(x_i-x_j) & \text{ if } i\neq j
\end{cases}
$$
Therefore, if $V''_{\{i,j\}}(u)\geq 0$ for any $u\in\dR$ and any
$i,j\in\{1,\ldots,n+1\}$, i.e. $V_{\{i,j\}}$ is convex, the
Gershgorin-Hadamard theorem implies that for any $x\in\dR^{n+1}$,
$\Hess{H}(x)\geq 0$ as a quadratic form, and thus $H$ is convex on
$\dR^{n+1}$.  Unfortunately, since $\sum_{j=1}^{n+1} \pd^2_{ij} H(x)=0$ for
any $i\in\{1,\ldots,n+1\}$, the null space of $\Hess{H}$ contains $\Vo_{n+1}$
and therefore, the measure $\mu$ on $\dR^{n+1}$ defined by
$d\mu(x):=\exp\PAR{-H(x)}\,dx$ cannot be normalised into a probability measure
since $Z_\mu:=\mu(\dR^{n+1})=+\infty$.  Nevertheless, suppose that there
exists a constant $\rho>0$ such that $V''_{\{i,j\}}(u)\geq \rho$ for any
$u\in\dR$ and any $i,j\in\{1,\ldots,n\}$.  Then, $u\in\dR\mapsto
V_{\{i,j\}}(u)-\rho u^2/2$ is convex and the latter implies that
$$
\Hess{H}(x) \geq \rho\,\bI_{n+1}-(n+1)^{-1}\rho\,\Mo_{n+1},
$$
as quadratic forms. Thus, by writing
$\dR^{n+1}=\dR\Vo_{n+1}\overset{\perp}{\oplus}\dH_n$ where $\dH_n$ is the
hyper-plane of equation $h_1+\cdots+h_{n+1}=0$, we get that the spectrum of
$\Hess{H}(x)$ is of the form
$$
\{0=\la_1(x)<\la_2(x)\leq\cdots\leq\la_{n+1}(x)\}
$$
with $\la_2(x)\geq n(n+1)^{-1}\rho$.  Hence, one can define the probability
measure $\si_M$ on $\dR^n$ as in \eqref{eq:defnum} for any $M$ in $\dR$.
Moreover $\si_M$ is uniformly log-concave with a constant $n(n+1)^{-1}\rho$ and
therefore the conditional measure $\mu_M$ can be defined from $\si_M$ as a
probability measure by equation \eqref{eq:mum2num}, despite the fact that
$\mu$ is not a probability measure on $\dR^{n+1}$.  The particular case
$V_{\{i,j\}}=V$ with $V$ even and uniformly convex is considered for example
in \cite{malrieu-projection}, in terms of the associated S.D.E., in order to
study the granular media equation.

As we have seen, when $H$ is uniformly strictly convex with a constant
$\rho>0$, the hypotheses of Proposition \ref{pr:num-mum} are full-filled and
hence, inequalities \eqref{in:ts-l2} and \eqref{in:ls-l2} hold. It is quite
natural to ask if \eqref{in:ts-l2} and \eqref{in:ls-l2} remains true for
symmetric but non convex Hamiltonians $H$.  In this direction, the
Bakry-\'Emery criterion allows the following perturbative statement due to
Ivan Gentil.  The proof, prototype of which can be found in
\cite{ledoux-lsrev}, is taken from \cite{hel-bod} and is postponed to section
\ref{se:pr:pr:ivan}.

\begin{prop}[Perturbative result]\label{pr:ivan}
  Let $H(x)=V(x_1)+\cdots+V(x_{n+1})$ with
  $$
  V(u)=\frac{u^2}{2}+F(u)
  $$
  where $F:\dR\to\dR$, and let $\si_M$ be the probability measure on
  $\dR^n$ defined by \eqref{eq:defnum}, namely
  $$
  \si_M(dx_1,\ldots,dx_n)
  =(Z_{\si_M})^{-1}
  \int_{\dR^{n}}\!\exp\BPAR{-\sum_{i=1}^n V(x_i)-V\BPAR{\MSX}}\,dx_1\cdots dx_n.
  $$
  Then, for \NI{F} small enough, there exists a positive constant $\cP$
  depending only on \NI{F} such that for any $n$, any $M$ and any smooth
  $f:\dR^n\to\dR$,
  \begin{equation}\label{in:ivan}
  \var{\si_M}{f} \leq \cP\,\moy{\si_M}{\ABS{\GR f}^2}.
  \end{equation}
\end{prop}
Proposition \ref{pr:ivan} remains valid if we replace, in the definition of
$\si_M$, the square function $u\mapsto u^2/2$ by a smooth convex function
$u\mapsto \Phi(u)$, provided that there exists real constants $\al$ and $\be$
such that $0<\al\leq\be\leq 2\al$ and $\al\leq \Phi''(u) \leq \be$ for every
$u\in\dR$.  The constant $\cP$ becomes in this case
$e^{2\,\OSC{F}}/\PAR{2\al\,e^{-2\,\OSC{F}}-\be}$ for
$\OSC{F}<\log\sqrt{2\al/\be}$.

The exchangeability of the underlying measure $\mu_M$ indicates that the
perturbative approach by mean of Helffer's method (cf.
\cite{helffer-98,helffer-99,helffer-99-2,hel-bod}) which sees $\si_M$ as a
quasi-product measure with small interactions is not relevant here: any
reduction of $F$ in the interaction term
$$
V\BPAR{\MSX}
$$
affects the product term $\sum_{i=1}^n V(x_i)$. Helffer's method was
essentially developed for spins systems with \emph{boundary conditions} for
which the measure is not exchangeable. For our measure $\si_M$, one can expect
in contrast that the symmetries of $H_M$ induce a stronger result, as for many
mean field models. In this direction, Landim, Panizo and Yau have recently
established in \cite{landim-panizo-yau} that $\mu_M$ satisfies inequalities
\eqref{in:ts-l2} and \eqref{in:ls-l2} when $H$ is of the form
$H(x)=V(x_1)+\cdots+V(x_{n+1})$ where $V(u)=u^2/2+F(u)$ with $F$ and $F'$
bounded and Lipschitz.  A simple example is given by $F(x)=P(\sin(Q(x)))$
where $P$ and $Q$ are fixed polynomials in $\dR[X]$. Their proof relies on
Lu-Yau's Markovian decomposition \cite{lu-yau} and on Local Central Limit
Theorem estimates \cite{kipnis-landim}.

Following closely \cite{landim-panizo-yau}, we are actually able to
show that measure $\si_M$ itself satisfies to \eqref{in:def-ts} and
\eqref{in:def-ls} with a constants which does not depend on $n$ and
$M$, as stated in our main result, which follows.

\begin{thm}\label{th:num-quad+bd}
  Let $H(x)=V(x_1)+\cdots+V(x_{n+1})$ with $V(u)=u^2/2+F(u)$ and let $\si_M$
  be the probability measure on $\dR^n$ defined by \eqref{eq:defnum}, namely
  $$
  \si_M(dx_1,\ldots,dx_n) =(Z_{\si_M})^{-1}
  \int_{\dR^{n}}\!\exp\BPAR{-\sum_{i=1}^n V(x_i)-V\BPAR{\MSX}}\,dx_1\cdots dx_n.
  $$
  Then, if $F$ is bounded and Lipschitz, there exists a positive constant
  $\cP$ depending only on \NI{F} and \NI{F'} such that for any $n$ and $M$ and
  any smooth $f:\dR^n\to\dR$,
  \begin{equation}\label{in:nu-quad+bd-ts}
  \var{\si_M}{f} \leq \cP\,\moy{\si_M}{\ABS{\GR f}^2}.
  \end{equation}  
  Moreover, if $F''$ is also bounded, there exists a positive constant
  $\cL$ depending only on \NI{F}, \NI{F'} and \NI{F''} such that for
  any $n$ and $M$ and any smooth $f:\dR^n\to\dR$,
  \begin{equation}\label{in:nu-quad+bd-ls}
  \ent{\si_M}{f^2} \leq \cL\, \moy{\si_M}{\ABS{\GR f}^2}.
  \end{equation}  
\end{thm}

As a Corollary, we recover from Proposition \ref{pr:num-mum} and
\eqref{in:compdirf} the $L^2$ factor for the Kawasaki dynamics (cf.
\eqref{in:ts-l2} and \eqref{in:ls-l2}) obtained by
\cite{landim-panizo-yau}.

The rest of the paper is divided as follows.  The first section gives
the proof of Proposition \ref{pr:ivan}, which relies only on the
Bakry-\'Emery criterion.  In Section \ref{se:prelims}, we give some
preliminaries to the proof of Theorem \ref{th:num-quad+bd}. Lemma
\ref{le:cov-estims} gives some covariance bounds taken from
\cite{landim-panizo-yau}. This Lemma allows us to derive the ``one
spin Lemma'' \ref{le:onespin} by a simple application of the
Bakry-\'Emery criterion.  Section \ref{se:qd+bd:ip} is devoted to the
derivation of the Poincar\'e inequality \eqref{in:nu-quad+bd-ts} and
section \ref{se:qd+bd:ls} to the derivation of the logarithmic Sobolev
counterpart \eqref{in:nu-quad+bd-ts}.  The proofs make heavy use of
the LCLT based estimates of \cite{landim-panizo-yau} throughout Lemmas
\ref{le:cvg} and \ref{le:ceg}, but our induction in $n$ is quite
different.

It is natural to ask if Theorem \ref{th:num-quad+bd} remains valid if the
quadratic potential $u^2/2$ is replaced by a uniformly strictly convex
potential $\Phi$. We believe that it is true. Recently, Caputo showed in
\cite{caputo2} that it is the case for the Poincar\'e inequality in Theorem
\ref{th:num-quad+bd}. His nice method makes crucial use of exchangeability,
but unfortunately, since it relies heavily on the spectral nature of
Poincar\'e's inequality, it does not give any clue to do the same for the
Logarithmic Sobolev inequality, and the second part of Theorem
\ref{th:num-quad+bd} remains thus inaccessible.  

In a sense, the exchangeability property plays a role similar to the one
played by mixing conditions in other models. Such exchangeable measures
``resemble'' to product ones, and this intuition is confirmed by a sort of
Kac's propagation of chaos since the finite dimensional marginals are close to
a product measure in high dimension, as we will see in Lemma
\ref{le:cov-estims}. Notice that in our exchangeable model with mean field
interaction, the covariance of any couple of spins decays linearly with the
total number of spins, whereas for spins systems with nearest neighbours
interaction and boundary conditions, the covariance decay holds exponentially.

The general study of Poincar\'e and Logarithmic Sobolev inequalities as in
Theorem \ref{th:num-quad+bd} for bounded ``diagonal'' perturbations of
non-exchangeable Hamiltonians is hard and remains an interesting open problem.
In an other direction, one can ask if our method remains valid for discrete
spins systems similar to those presented in \cite{martinelli-99}. It is not
clear at all for us. Finally, we believe that concentration of measure
inequalities can help to simplify the derivation of large deviations like
estimates in \cite{landim-panizo-yau} necessary to derive the Logarithmic
Sobolev inequality.

\section{Proof of Proposition \ref{pr:ivan}}\label{se:pr:pr:ivan}

We give here a proof of Proposition \ref{pr:ivan}, which relies only on the
Bakry-\'Emery criterion. Let $\si_M^*$ the probability measure on $\dR^n$
defined by
$$
(Z_{\si_M^*})^{-1}\,\exp\BPAR{-\sum_{i=1}^n
  V(x)-\frac{1}{2}\,\BPAR{\MSX}^2\,}\,dx_1\cdots dx_n.
$$
If $\si_M^*$ satisfies a Poincar\'e inequality with a constant $c>0$,
then $\si_M$ satisfies a Poincar\'e inequality with a constant
$c\exp\PAR{2\,\OSC{F}}$.  Now, for any smooth function
$f:\dR^n\to\dR$,
$$
\moy{\si_M^*}{(\GI f)^2}=\sum_{i,j=1}^n\moy{\si_M^*}{\ABS{\pd^2_{ij} f}^2}
+\bmoy{\si_M^*}{\sum_{i=1}^n\PAR{1+F''(x_i)}\,\ABS{\pd_i f}^2}
+\bmoy{\si_M^*}{\BPAR{\sum_{i=1}^n \pd_i f}^2}.
$$
In the other hand, for any $i\in\{1,\ldots,n\}$ and any
$x_1,\ldots,x_{i-1},x_{i+1},\ldots,x_n$, the Bakry-\'Emery criterion gives
that the one dimensional probability measure
$$
\rho_i(dx_i):=(Z_{\rho_i})^{-1}\,
\exp\BPAR{-V(x_i)-\frac{1}{2}\BPAR{\MSX}^2\,}\,dx_i
$$
satisfies a Poincar\'e inequality with a constant $(1/2)\,\exp(2\,\OSC{F})$,
hence, by the Bakry-\'Emery criterion applied reversely, we get for any smooth
function $f:\dR^n\to\dR$, by summing over $i$
$$
\sum_{i=1}^n\moy{\rho_i}{\ABS{\pd^2_{ii} f}^2} +
\sum_{i=1}^n\moy{\rho_i}{\PAR{2+F''(x_i)}\,\ABS{\pd_i f}^2} \geq
2\,e^{-2\,\OSC{F}}\sum_{i=1}^n\moy{\rho_i}{\ABS{\pd_i f}^2}.
$$
Notice that
$\rho_i=\LAW{\si_M^*}{x_i\,\vert\,x_1,\ldots,x_{i-1},x_{i+1},\ldots,x_n}$.
Therefore, by taking the expectation with respect to $\si_M^*$, we get
\begin{align*}
  \moy{\si_M^*}{(\GI f)^2} 
  &\geq \PAR{2\,e^{-2\,\OSC{F}}-1}\,\sum_{i=1}^n\moy{\si_M^*}{\ABS{\pd_i f}^2}\\
  &=: \PAR{2\,e^{-2\,\OSC{F}}-1}\,\moy{\si_M^*}{\ABS{\GR f}^2}.
\end{align*}
Thus, for $\OSC{F}$ sufficiently small ($<\log\sqrt{2}$), one can take
$$
\cP=\frac{e^{2\,\OSC{F}}}{2\,e^{-2\,\OSC{F}}-1},
$$
which is optimal when $F\equiv0$ (pure Gaussian case).

\section{Preliminaries to the proof of Theorem \ref{th:num-quad+bd}}
\label{se:prelims}

Let $\ga_{n,M}$ the Gaussian measure of mean $M/(n+1)$ and covariance matrix
$\PAR{\bI_n+\Mo_n}^{-1}$. If $B(x):=\sum_{i=1}^n
F(x_i)+F\PAR{M-x_1-\cdots-x_n}$, one can write
$$
d\si_M(x_1,\ldots,x_n)
=(Z_{n,M})^{-1}\,\exp\PAR{-B(x)}\,d\ga_{n,M}(x_1,\ldots,x_n).
$$
Thus, $\si_M$ is a bounded perturbation of $\ga_{n,M}$, which is
log-concave with a constant $\rho$ equal to $1$, and therefore,
$\si_M$ satisfies to Poincar\'e and logarithmic Sobolev inequalities
with constants depending only on \NI{B} (i.e. $\NI{F}$ and $n$). Our
goal is to show that the dependence in $n$ can be dropped by taking
into account \NI{F'} and \NI{F''}. The presence of the bounded part
$F$ in $V$ and the non-product nature of $\si_M$ does not allow any
direct approach based on the Bakry-\'Emery criterion.

Observe that $\cov{\ga_{n,M}}{x_1}{x_2}=-(n+1)^{-1}$, and we can then expect
the same decrease in $n$ for $\cov{\si_M}{V'(x_1)}{V'(x_2)}$. This is actually
the case, as stated in the following Lemma. Notice that since $\si_M$ is
exchangeable and since $\MSX$ and $x_i$ have the same law under $\si_M$, we
have $\var{\si_M}{x_1}=-n\cov{\si_M}{x_1}{x_2}$, as for $\ga_{n,M}$.

\begin{lem}\label{le:cov-estims}
  Let $\si_M$ be the probability measure on $\dR^n$ ($n>2$) defined in Theorem
  \ref{th:num-quad+bd} and $\mu_M$ the associated conditional measure defined
  by \eqref{eq:defmum}. Assume that $F$ and $F'$ are bounded, then there
  exists a constant $C>0$ depending only on \NI{F} and \NI{F'} such that for
  any $M\in\dR$
  \begin{equation}
    \label{in:cov-estim1}
    \ABS{\cov{\si_M}{V'(x_1)}{V'(x_2)}} =
    \ABS{\cov{\mu_M}{V'(x_1)}{V'(x_2)}} \leq \frac{C}{n},
  \end{equation}
  and
  \begin{equation}\label{in:cov-estim2}
  \begin{split}
    \var{\si_M}{\sum_{i=1}^n V'(x_i)+V'(\MSX)}
    &=\var{\si_M}{\sum_{i=1}^n F'(x_i)+F'(\MSX)}\\
    &=\var{\mu_M}{\sum_{i=1}^{n+1} F'(x_i)} \leq nC.
  \end{split}
  \end{equation}
\end{lem}

\begin{proof}
  Inequality \eqref{in:cov-estim2} follows from \eqref{eq:mum2num} and
  \cite[Corollary 5.4]{landim-panizo-yau}.
  For \eqref{in:cov-estim1}, just write
  $$
  \cov{\si_M}{V'(x_1)}{V'(x_2)}
  =\cov{\si_M}{x_1}{x_2}
  +2\,\cov{\si_M}{x_1}{F'(x_2)}
  +\cov{\si_M}{F'(x_1)}{F'(x_2)},
  $$
  and use \eqref{eq:mum2num} and \cite[Corollary 5.3]{landim-panizo-yau} to
  estimate each term.  Actually, one can derive the estimates of
  $\cov{\si_M}{x_1}{F'(x_2)}$ and $\cov{\si_M}{x_1}{x_2}$ directly by using
  the symmetries of $\si_M$.
\end{proof}

Inequality \eqref{in:cov-estim1} of Lemma \ref{le:cov-estims} allows
us to establish the following one spin result, which is the first step
in our proof of Poincar\'e and logarithmic Sobolev inequalities for
$\si_M$ by induction on $n$ by mean of the Lu-Yau Markovian
decomposition. In the other hand, inequality \eqref{in:cov-estim2}
will be usefull, as we will see in sections \ref{se:qd+bd:ip} and
\ref{se:qd+bd:ls}, for the induction itself.

\begin{lem}[One spin Lemma]\label{le:onespin}
  Let $\si_M$ be the probability measure on $\dR^n$ defined in Theorem
  \ref{th:num-quad+bd}. If $F$ is bounded and Lipschitz, there exists a
  constant $A>0$ depending only on \NI{F} and \NI{F'} and not on $n$ and $M$
  such that for any $n$ and $M$ and any smooth $f:\dR\to\dR$,
  $$
  \ent{\si_M}{f(x_1)^2} \leq 2A\,\moy{\si_M}{f'(x_1)^2},
  $$
  and
  $$
  \var{\si_M}{f(x_1)} \leq A\,\moy{\si_M}{f'(x_1)^2}.
  $$
\end{lem}

\begin{proof}[Proof of Lemma \ref{le:onespin}]\label{le:pr:onespin}
  As we already noticed, it is clear that the desired inequalities are true
  with a constant depending on $n$ and $\NI{F}$, so we just have to see what
  happens for large values of $n$. We have in mind the use of the Bakry-\'Emery
  criterion. The Hamiltonian of the probability measure in $x_1$ is given by
  \begin{multline*}
    \vphi_{M,n}(x_1):=V(x_1)+\log Z_{M,n} 
    - \log \int\!\exp\PAR{-\sum_{i=2}^n V(x_i)-V\BPAR{\MSX}}\, dx_2\cdots dx_n.
  \end{multline*}
  We first observe that we can forget the $F(x_1)$ part in $V(x_1)$, which is
  payed by a factor $\exp(2\,\OSC{F})$ in $A$. Hence, we simply have, after an
  integration by parts
  \begin{equation*}
    \vphi_{M,n}''(x_1)=1-\cov{\si_{M-x_1}(dx_2,\ldots,dx_n)}{V'(x_2)}{V'(x_3)}.
  \end{equation*}
  Now, \eqref{in:cov-estim1} gives $\vphi_{M,n}''(x_1) \geq 1-Cn^{-1}$, where
  $C$ is a positive constant depending only on \NI{F} and \NI{F'} and not on
  $n$ and $M$. Thus, we are able to apply the Bakry-\'Emery criterion for large
  values of $n$. Hence, the proof is completed, with a constant $A$ depending
  only on \NI{F} and \NI{F'} and not on $M$ and $n$.
\end{proof}

Obviously, one can replace $x_1$ in $f$ and $f'$ by $M-x_1-\cdots-x_n$ or by
any $x_i$ for $i\in\{1,\ldots,n\}$. Moreover, according to
\eqref{eq:mumeqnum}, one can replace $\moyf{\si_M}$ by $\moyf{\mu_M}$.

\section{Derivation of the Poincar\'e inequality}\label{se:qd+bd:ip}

This section is devoted to the derivation of inequality
\eqref{in:nu-quad+bd-ts} of Theorem \ref{th:num-quad+bd}. The proof
relies on the one spin Lemma \ref{le:onespin} and on the crucial Lemma
\ref{le:cvg} which allows us to use the Lu-Yau Markovian
decomposition.

\begin{proof}[Proof of \eqref{in:nu-quad+bd-ts}]\label{th:pr:nu-quad+bd:ts} 
  As we already noticed, the result is true with a constant depending
  on $n$, so that if we denote by $\cP_n$ the maximum of best Poincar\'e
  constants in dimension less than or equal to $n$, we just have to
  show that the non decreasing sequence of constants $(\cP_n)_{n\geq
    1}$ is bounded.
  
  Let us denote by $\si$ the measure $\si_M$ and by $\sik{k}$ the measure
  $\si_M$ given $x_1,\ldots,x_k$ for $k\in\{0,\ldots,n\}$ and by $f_k$ the
  conditional expectation
  $$
  \moy{\si}{f\vert x_1,\ldots,x_k}=\moy{\sik{k}}{f}.
  $$
  Notice that $\sik{k}$ is nothing else but
  $\si_{M-x_1-\cdots-x_k}(dx_{k+1},\ldots,dx_n)$. Moreover, $f_n=f$ and by
  convention $\si^{(0)}:=\si$ and thus $f_0=\moy{\mu}{f}$.  For a \emph{fixed}
  function $f$, we can always choose the order of the coordinates
  $x_1,\ldots,x_n$ such that $\moy{\si}{\ABS{\pd_k f}^2}$ becomes a \emph{non
    increasing} sequence in $k\in\{1,\ldots,n\}$. This gives
  $$
  \sum_{i=k+1}^n \frac{1}{n-k}\,\moy{\si}{\ABS{\pd_i f}^2} \leq
  \moy{\si}{\ABS{\pd_{k+1} f}^2}.
  $$
  Following Lu-Yau \cite{lu-yau}, we have the following Markovian
  decomposition of the variance
  $$
  \var{\si}{f} := \moy{\si}{f^2}-\moy{\si}{f}^2
  =\sum_{k=1}^n\moy{\si}{(f_k)^2-(f_{k-1})^2}
  =\sum_{k=1}^n\moy{\si}{\var{\sik{k-1}}{f_k}}.
  $$
  Since measure $\sik{k-1}$ integrates coordinates $x_k,\ldots,x_n$ and
  function $f_k$ depends only on coordinates $x_1,\ldots,x_k$, the quantity
  $\var{\sik{k-1}}{f_k}$ is actually a variance for a one spin function.
  Therefore, by the one spin Lemma \ref{le:onespin}, there exists a constant
  $A>0$ depending on \NI{F} and \NI{F'} but not on $n$ and $M$ such that
  $$
  \var{\si}{f} \leq A\,\sum_{k=1}^n \moy{\si}{\ABS{\pd_k f_k}^2}.
  $$
  Our aim is to express the right hand side of the previous inequality in
  terms of $\ABS{\pd_k f}^2$. Notice that the $k=n$ term in the sum is trivial
  since $f_n=f$. By definition of $f_k$, we get for any $k\in\{1,\ldots,n-1\}$
  $$
  \pd_k f_k=\moy{\sik{k}}{\pd_k f}-\bcov{\sik{k}}{f}{V'(\MSX)}.
  $$
  At this stage, we notice that by $n-k$ integrations by parts, we have
  $$
  \bcov{\sik{k}}{f}{V'(\MSX)}
  =\frac{1}{n-k}\, \sum_{i=k+1}^n\cov{\sik{k}}{f}{V'(x_i)}
  -\frac{1}{n-k}\,\sum_{i=k+1}^n \moy{\sik{k}}{\pd_i f}.
  $$
  Therefore, we can write by denoting $S_k:=\sum_{i=k+1}^n V'(x_i)+V'(\MSX)$
  $$
  \pd_k f_k =\moy{\sik{k}}{\pd_k f} -\frac{1}{n-k+1}\,\cov{\sik{k}}{f}{S_k}
  +\frac{1}{n-k+1}\,\sum_{i=k+1}^n \moy{\sik{k}}{\pd_i f}.
  $$
  Now, by the Cauchy-Schwarz inequality
  $$
  \ABS{\pd_k f_k}^2 
  \leq 3\,\moy{\sik{k}}{\ABS{\pd_k f}^2}
  +\frac{3}{(n-k)^2}\,\cov{\sik{k}}{f}{S_k}^2
  +\frac{3}{n-k}\,\sum_{i=k+1}^n\moy{\sik{k}}{\ABS{\pd_i f}^2}.
  $$
  This gives by summing over all $k$ in $\{1,\ldots,n-1\}$ (the case $k=n$
  is trivial)
  \begin{align*}
    \sum_{k=1}^{n-1}\moy{\si}{\ABS{\pd_k f_k}^2}
    \leq\;&3\,\moy{\si}{\ABS{\GR f}^2}
          +3\,\sum_{k=1}^{n-1}
          \frac{1}{(n-k)^2}\,\moy{\si}{\cov{\sik{k}}{f}{S_k}^2}\\
          &+3\,\sum_{k=1}^{n-1}
          \frac{1}{n-k}\,\sum_{i=k+1}^n\moy{\si}{\ABS{\pd_if}^2}.
  \end{align*}
  The monotonicity of $\moy{\si}{\ABS{\pd_i f}^2}$ yields
  $$
  \sum_{k=1}^{n-1}
  \moy{\si}{\ABS{\pd_k f_k}^2} \leq\;6\,\moy{\si}{\ABS{\GR f}^2}
  +3\,\sum_{k=1}^{n-1}\frac{1}{(n-k)^2}\,\moy{\si}{\cov{\sik{k}}{f}{S_k}^2}.
  $$
  By inequality \eqref{in:cvg} of Lemma \ref{le:cvg}, there exists a
  positive constant $C$ depending only on \NI{F} and \NI{F'} such that for any
  $\veps>0$, there exists a positive constant $C_\veps$ depending only on
  \NI{F}, \NI{F'} and $\veps$ such that for any $k\in\BRA{1,\ldots,n-1}$
  $$
  \cov{\sik{k}}{f}{S_k}^2 \leq (C_\veps+\veps (n-k)C)\,\var{\sik{k}}{f}
  +(n-k)C_\veps\,\sum_{i=k+1}^n\moy{\sik{k}}{\ABS{\pd_i f}^2}.
  $$
  Therefore, by the monotonicity of $\moy{\si}{\ABS{\pd_i f}^2}$
  again
  \begin{equation}\label{in:pr:estim:grad-grad-var-var}
  \begin{split}
    \sum_{k=1}^{n-1}\moy{\si}{\ABS{\pd_k f_k}^2}
    \leq\;&C'_\veps\,\moy{\si}{\ABS{\GR f}^2}
    +C'_\veps\,\sum_{k=1}^{n-1}%
    \frac{1}{(n-k)^2}\,\moy{\si}{\var{\sik{k}}{f}} \\
    &+\veps
    C'\,\sum_{k=1}^{n-1}\frac{1}{n-k}\,\moy{\si}{\var{\sik{k}}{f}}.
  \end{split}
  \end{equation}
  Recall that $\cP_n$ is the maximum of best Poincar\'e constants in dimension
  less than or equal to $n$. The last sum of the right hand side (RHS) of
  \eqref{in:pr:estim:grad-grad-var-var} can be bounded above as follows
  $$
  \sum_{k=1}^{n-1}\frac{1}{n-k}\moy{\si}{\var{\sik{k}}{f}} 
  \leq \cP_{n-1} \moy{\si}{\ABS{\GR f}^2}.
  $$
  It remains to examine the first sum of the RHS of
  \eqref{in:pr:estim:grad-grad-var-var}.  The Jensen inequality yields
  $$
  \moy{\si}{\var{\sik{k}}{f}} \leq \var{\si}{f},
  $$
  and therefore, we get for any $p\in\{1,\ldots,n-1\}$
  \begin{align*}
    \sum_{k=1}^{n-1}\frac{1}{(n-k)^2}\,\moy{\si}{\var{\sik{k}}{f}}
    &=\var{\si}{f}\sum_{k=1}^{n-p-1}\frac{1}{(n-k)^2}
    +\sum_{k=n-p}^{n-1}\frac{1}{(n-k)^2}\,\moy{\si}{\var{\sik{k}}{f}}\\
    &\leq\var{\si}{f}\sum_{k=p+1}^{n-1}\frac{1}{k^2}
    +\sum_{k=1}^{p}\frac{1}{k^2}\,\moy{\si}{\var{\sik{n-k}}{f}}.
  \end{align*}
  At this stage, we observe that for every $k$ in $\{1,\ldots,p\}$,
  $$
  \moy{\si}{\var{\sik{n-k}}{f}} 
  \leq \cP_{p} \sum_{i=n-k+1}^n \moy{\si}{\ABS{\pd_i f}^2}
  \leq p\cP_{p}\,\moy{\si}{\ABS{\pd_{n-p+1} f}^2}.
  $$
  We are now able to collect our estimates of the RHS of
  \eqref{in:pr:estim:grad-grad-var-var}. Putting all together, we have
  obtained that
  $$
  \sum_{k=1}^{n-1}\moy{\si}{\ABS{\pd_k f_k}^2}
  \leq 
  (C'_\veps+p\pi^2\cP_{p}C'_\veps+\veps C'\cP_{n-1})\,\moy{\si}{\ABS{\GR f}^2}
  +(C'_\veps R_p)\var{\si}{f},
  $$
  where $R_p:=\sum_{k=p+1}^{n-1}k^{-2}$. Therefore, for some
  $C''_{p,\veps}>0$,
  $$
  (1-AC'_\veps R_p)\,\var{\si}{f}
  \leq (C''_{p,\veps}+\veps AC'\cP_{n-1})\,\moy{\si}{\ABS{\GR f}^2}.
  $$
  Now, we may choose $\veps < 1/(AC')$ and then $p$ large enough (always
  possible when $n$ is sufficiently large) to ensure that
  $$
  R_p < \min\PAR{\frac{1}{AC'_\veps},\frac{1-\veps AC'}{AC'_\veps}}.
  $$
  This gives two positive constants $\al$ and $\be$ with $\be<1$ depending
  only on \NI{F} and \NI{F'} such that for large values of $n$, one has
  $\cP_n \leq \al + \be\,\cP_{n-1}$, and therefore $\sup_{n} \cP_n < +\infty$.
\end{proof}

Let us give now the crucial Lemma which allows us to use the Markovian
decomposition of Lu-Yau, by splitting the covariance term into a
variance term and a gradient term. The proof makes heavy use of
estimates taken from \cite{landim-panizo-yau}.

\begin{lem}\label{le:cvg} 
  Let $\si_M$ be the probability measure on $\dR^n$ defined in Theorem
  \ref{th:num-quad+bd}. Assume that $F$ is bounded and Lipschitz, then there
  exists a positive constant $C$ depending only on $\NI{F}$ and $\NI{F'}$ such
  that for any $\veps>0$, there exists a positive constant $C_\veps$ depending
  only on $\NI{F}$, $\NI{F'}$ and $\veps$ such that for any $n\in\dN^*$, any
  $M\in\dR$ and any smooth function $f:\dR^n\to\dR$
 \begin{equation}\label{in:cvg}   
  \cov{\si_M}{f}{S}^2
  \leq (C_\veps+\veps n C)\,\var{\si_M}{f}
  +nC_\veps\,\moy{\si_M}{\ABS{\GR f}^2},
 \end{equation}
 where $S:=\sum_{i=1}^n V'(x_i)+V'(\MSX)$.
\end{lem}

\begin{proof}[Proof of Lemma \ref{le:cvg}]\label{le:pr:cvg}
  Notice that we just have to study what happens for small values of $\veps$
  and large values of $n$, since for any $\veps>0$ and any $n\leq n_\veps$, we
  get by the Cauchy-Schwarz inequality and \eqref{in:cov-estim2} that
  $$
  \cov{\si_M}{f}{S}^2 \leq n_\veps C\,\var{\si_M}{f} =:
  C_\veps\,\var{\si_M}{f}.
  $$
  We have in mind the use of the partitioning result of
  \cite{landim-panizo-yau}. If ${\mu_M}$ denotes the conditional measure on
  $\dR^{n+1}$ associated to $\si_M$ as in \eqref{eq:mum2num}, we have
  $$
  \cov{{\si_M}}{f}{S(x_1,\ldots,x_n)}^2
  =\bcov{{\mu_M}}{f}{\sum_{i=1}^{n+1}F'(x_i)}^2.
  $$
  Now, for $n$ large enough, one can then subdivide the set
  $\{1,\ldots,n+1\}$ into $\ell$ adjacent subsets $I_i$ of size $K$ or $K+1$.
  We have in mind to take $K^{-1}\leq\veps$, which is always possible when $n$
  is large enough. We can write with this decomposition
  $$
  \bcov{{\mu_M}}{f}{\sum_{i=1}^{n+1} F'(x_i)} =
  \bcov{{\mu_M}}{f}{\sum_{i=1}^\ell\sum_{k\in I_i} F'(x_k)}.
  $$
  For any $(i,x)\in\BRA{1,\ldots,\ell}\times\dR^{n+1}$, we define the
  ``total spin on $I_i$'' by $M_i(x):=\sum_{k\in I_i} x_k$.  On $\dR^{I_i}$,
  one can define the conditional measure $\mu_{M_i}$ with total spin $M_i$,
  as in \eqref{eq:defmum}. To lighten the notations, we denote this measure by
  $\muMi$.  We get from the latter by the Cauchy-Schwarz inequality
  \begin{equation}\label{in:cvg-step-1}
    \begin{split}
      \bcov{{\mu_M}}
      {f}{\sum_{i=1}^{n+1}F'(x_i)}^2 
      \leq\;&2\,\bcov{{\mu_M}}{f}
      {\sum_{i=1}^\ell\sum_{k\in I_i}\PAR{F'(x_k)-\moy{\muMi}{F'}}}^2\\
      &+ 2\,\bcov{{\mu_M}} {f}{\sum_{i=1}^\ell\ABS{I_i}\moy{\muMi}{F'}}^2
    \end{split}
  \end{equation}
  By the Cauchy-Schwarz inequality again and by \eqref{eq:mumeqnum},
  the second term of the RHS of \eqref{in:cvg-step-1} can be bounded
  above by
  $$
  \var{\si_M}{f}\bvar{\mu_M}{\sum_{i=1}^\ell\ABS{I_i}\moy{\muMi}{F'}}.
  $$
  Now, according to \cite[ineq. (3.10)]{landim-panizo-yau}, the last
  variance in the RHS is bounded above by $nC/K$ for $n$ sufficiently large,
  which can be rewritten as $\veps nC$.  We turn now to the control of the
  first term of the RHS of \eqref{in:cvg-step-1}. Since
  $\moyf{{\mu_M}}=\moyf{{\mu_M}} \circ \moyf{\muMi}$, we get
  $$
  \bcov{{\mu_M}}
  {f}
  {\sum_{i=1}^\ell\sum_{k\in I_i}\PAR{F'(x_k)-\moy{\muMi}{F'}}}
  =\sum_{i=1}^\ell\bmoy{{\mu_M}}{\bcov{\muMi}{f}{\sum_{k\in I_i} F'(x_k)}}.
  $$
  Thus, the Cauchy-Schwarz inequality yields
  $$
  \bcov{{\mu_M}}
  {f}
  {\sum_{i=1}^\ell\sum_{k\in I_i}\PAR{F'(x_k)-\moy{\muMi}{F'}}}^2
  \leq
  \ell\,
  \sum_{i=1}^\ell\bmoy{{\mu_M}}{\bcov{\muMi}{f}{\sum_{k\in I_i} F'(x_k)}^2}.
  $$
  Again by the Cauchy-Schwarz inequality, we get
  $$
  \bcov{\muMi}{f}{\sum_{k\in I_i} F'(x_k)}^2 
  \leq
  \var{\muMi}{f}\bvar{\muMi}{\sum_{k\in I_i}F'(x_k)}.
  $$
  By virtue of \eqref{in:cov-estim2} applied to $\muMi$, we obtain
  $$
  \bcov{\muMi}{f}{\sum_{k\in I_i} F'(x_k)}^2
  \leq C\,\ABS{I_i}\,\var{\muMi}{f}.
  $$
  Now, for any $i$, let $r_i=\max\{k,\,k\in I_i\}$
  and $J_i:=I_i\bs\{r_i\}$ and $\siMi$
  the probability measure on $\dR^{J_i}$ associated with the Hamiltonian
  $$
  \sum_{k\in J_i} V(x_k)+V(M_i-\sum_{k\in J_i} x_k).
  $$
  Equation \eqref{eq:mumeqnum} simply gives
  $$
  \var{\muMi}{f} = \var{\siMi}{f(\vphi_i(x))}, 
  $$
  where $\vphi_i:\dR^{n+1}\to\dR^n$ is defined by
  $$
  (\vphi_i(x))_k
  :=
  \begin{cases}
    x_k & \text{if $k\neq r_i$}\\
    M_i-\sum_{l\in J_i} x_l & \text{if $k=r_i$}
  \end{cases}
  $$
  Recall that $\cP_K$ is the maximum of the best Poincar\'e constants for
  $\si_M$ in dimensions less than or equal to $K$. We get by definition of
  $\cP_K$ that
  \begin{align*}
    \var{\siMi}{f} 
    &\leq\cP_K\sum_{k\in J_i}
    \moy{\siMi}{\ABS{(\pd_k f)(\vphi_i)-(\pd_{r_i}f)(\vphi_i)}^2}\\
    &=\cP_K\sum_{k\in J_i}
    \moy{\muMi}{\ABS{\pd_k f-\pd_{r_i}f}^2}.
  \end{align*}
  Hence, by the Cauchy-Schwarz inequality, we get
  $$
  \var{\siMi}{f}
  \leq
  2\cP_K\,\bmoy{\muMi}{\sum_{k\in J_i} \ABS{\pd_k f}^2} +
  2(\ABS{I_i}-1)\cP_K\,\moy{\muMi}{\ABS{\pd_{r_i} f}^2},
  $$
  Summarising, since $\cP_K$ depends only on $K$, \NI{F}, \NI{F'},
  we have obtained that the first term of the right hand side of
  \eqref{in:cvg-step-1} is bounded above by
  $$
  nC_K\cP_K\,\moy{\mu_M}{\ABS{\GR f}^2},
  $$
  which can be rewritten by virtue of \eqref{eq:mumeqnum} as
  $n\,C'_\veps\,\moy{{\si_M}}{\ABS{\GR f}^2}$. This concludes the proof of
  \eqref{in:cvg} and Lemma \ref{le:cvg}.
\end{proof}

\section{Derivation of the Logarithmic Sobolev inequality}\label{se:qd+bd:ls}

This section is devoted to the derivation of inequality
\eqref{in:nu-quad+bd-ls} of Theorem \ref{th:num-quad+bd}. As for the
Poincar\'e inequality \eqref{in:nu-quad+bd-ts}, the proof relies on the
one spin Lemma \ref{le:onespin} and on a crucial Lemma \ref{le:ceg}
which allows us to use the Lu-Yau Markovian decomposition.

\begin{proof}[Proof of the logarithmic Sobolev inequality
  \eqref{in:nu-quad+bd-ls} of Theorem
  \ref{th:num-quad+bd}]\label{th:pr:nu-quad+bd:ls}
  
  We follow here the same scheme used for the Poincar\'e inequality.  For any
  smooth non negative function $g:\dR^n\to\dR^+$, we have the
  following decomposition of the entropy
  \begin{align*}
    \ent{\si}{g} :&= \moy{\si}{g\log g}-\moy{\si}{g}\log\moy{\si}{g}\\
    &=\sum_{k=1}^n \moy{\si}{g_k\log g_k-g_{k-1}\log g_{k-1}}\\
    &=\sum_{k=1}^n \moy{\si}{\ent{\sik{k-1}}{g_k}}.
  \end{align*}
  Alike for the variance, measure $\sik{k-1}$ integrates on
  $x_k,\ldots,x_n$ and function $f_k$ depends only on
  $x_1,\ldots,x_k$, so that $\ent{\sik{k-1}}{g_k}$ is actually an
  entropy for a one spin function.  Therefore, by the one spin Lemma
  \ref{le:onespin}, there exists a positive constant $A$ depending on
  \NI{F} and \NI{F'} but not on $n$ and $M$ such that
  $$
  \ent{\si}{g}
  \leq 2A\,\sum_{k=1}^n \bmoy{\si}{\frac{\ABS{\pd_k g_k}^2}{4g_k}},
  $$
  By taking $g=f^2$ for a smooth function $f:\dR^n\to\dR$, we get
  $$
  \ent{\si}{f^2}
  \leq 2A\,\sum_{k=1}^n \bmoy{\si}{\frac{\ABS{\pd_k (f^2)_k}^2}{4(f^2)_k}}.
  $$
  By imitating the method used for the Poincar\'e inequality, we get
  that
  $$
  \frac{\ABS{\pd_k (f^2)_k}^2}{4(f^2)_k}
  \leq 3\,\frac{\ABS{\moy{\sik{k}}{f\pd_k f}}^2}{\moy{\sik{k}}{f^2}}
  +\frac{6}{(n-k)^2}\,\frac{\cov{\sik{k}}{f^2}{S_k}^2}{\moy{\sik{k}}{f^2}}
  +\frac{3}{n-k}\,\sum_{i=k+1}^n 
  \frac{\ABS{\moy{\sik{k}}{f\pd_i f}}^2}{\moy{\sik{k}}{f^2}}.
  $$
  The Cauchy-Schwarz inequality yields
  $$
  \frac{\ABS{\moy{\sik{k}}{f\pd_k f}}^2}{\moy{\sik{k}}{f^2}}
  \leq\moy{\sik{k}}{\ABS{\pd_k f}^2}.
  $$
  Therefore, the Jensen inequality and the monotonicity of
  $\moy{\si}{\ABS{\pd_i f}^2}$ yield
  $$
  \sum_{k=1}^{n-1}\moy{\si}{\frac{\ABS{\pd_k (f^2)_k}^2}{4(f^2)_k}}
  \leq\; 6\,\moy{\si}{\ABS{\GR f}^2}
  +6\,\sum_{k=1}^{n-1}\frac{1}{(n-k)^2}\,
  \moy{\si}{\frac{\cov{\sik{k}}{f^2}{S_k}^2}{\moy{\sik{k}}{f^2}}}.
  $$
  By inequality \eqref{in:ceg} of Lemma \eqref{le:ceg}, there
  exists a positive constant $C$ depending only on \NI{F}, \NI{F'} and
  \NI{F''} such that for any $\veps>0$, there exists a positive
  constant $C_\veps$ depending only on \NI{F}, \NI{F'}, \NI{F''} and
  $\veps$ such that for any $n$ and $M$
  $$
  \frac{\cov{\sik{k}}{f^2}{S_k}^2}{\moy{\sik{k}}{f^2}}
  \leq (C_\veps+\veps(n-k)C)\,\ent{\sik{k}}{f^2}
  +(n-k)C_\veps\,\sum_{i=k+1}^n\moy{\sik{k}}{\ABS{\pd_i f}^2}.
  $$
  Hence, we are now able to proceed as the same way as for the Poincar\'e
  inequality.
\end{proof}

As for the derivation of the Poincar\'e inequality, we give now the
crucial Lemma which allows us to use the Markovian decomposition of
Lu-Yau.

\begin{lem}\label{le:ceg} 
  Let $\si_M$ be the probability measure on $\dR^n$ defined in Theorem
  \ref{th:num-quad+bd}. Assume that $F$, $F'$ and $F''$ are bounded, then
  there exists a positive constant $C$ depending only on $\NI{F}$, $\NI{F'}$
  and $\NI{F''}$ such that for any $\veps>0$, there exists a positive constant
  $C_\veps$ depending only on $\NI{F}$, $\NI{F'}$, $\NI{F''}$ and $\veps$ such
  that for any $n\in\dN^*$, any $M\in\dR$ and any smooth function
  $f:\dR^n\to\dR$ such that $\moy{\si_M}{f^2}=1$
  \begin{equation}\label{in:ceg}   
    \cov{\si_M}{f^2}{S}^2
    \leq (C_\veps+\veps n C)\,\ent{\si_M}{f^2}
    +nC_\veps\,\moy{\si_M}{\ABS{\GR f}^2},
  \end{equation}
  where $S(x):=\sum_{i=1}^n V'(x_i)+V'(\MSX)$.
\end{lem}

\begin{proof}[Proof of Lemma \ref{le:ceg}]\label{le:pr:ceg}
  We follow the same scheme as for \eqref{in:cvg}, by replacing the
  Cauchy-Schwarz inequality by the entropy inequality.  Since $f^2$ is
  a density with respect to $\si_M$, we can write
  $$
  \cov{\si_M}{f^2}{S}=\moy{\si_M}{\PAR{S-\moy{\si_M}{S}}\,f^2},
  $$
  and hence, we get by the entropy inequality that for any $\be>0$
  $$
  \cov{\si_M}{f^2}{S} \leq
  \be^{-1}\,\log\moy{\si_M}{\exp\PAR{\be\PAR{S-\moy{\si_M}{S}}}}
  +\be^{-1}\,\moy{\si_M}{f^2\log f^2}.
  $$
  By \eqref{eq:mum2num} and \cite[Lemma 6.1]{landim-panizo-yau},
  the first term of the right hand side is bounded above by $nC\be$
  where $C$ depends only on \NI{F}, and \NI{F''}. This yields by
  considering the minimum in $\be>0$
  $$
  \cov{\si_M}{f^2}{S}^2\leq nC\,\moy{\si_M}{f^2\log f^2}.
  $$
  Thus, for any fixed $\veps>0$, we just have to study what happens
  for large values of $n$ since $nC\leq n_\veps C=:C_\veps$ for $n\leq
  n_\veps$.  After rewriting \eqref{in:ceg} in terms of ${\mu_M}$, we
  get by Cauchy-Schwarz's inequality
  \begin{equation}\label{in:ceg-step-1}
    \begin{split}
    \bcov{{\mu_M}}
    {f^2}{\sum_{i=1}^{n+1}F'(x_i)}^2 
    \leq\;&
    2\,\bcov{{\mu_M}}
    {f^2}{\sum_{i=1}^\ell\sum_{k\in I_i}\PAR{F'(x_k)-\moy{\muMi}{F'}}}^2\\
    &+2\,\bcov{{\mu_M}}{f^2}{\sum_{i=1}^\ell\ABS{I_i}\moy{\muMi}{F'}}^2.
    \end{split}
  \end{equation}
  Let us treat the first term of the right hand side of
  \eqref{in:ceg-step-1}. It can be rewritten as
  $$
  2\,\sum_{i=1}^\ell
  \bmoy{\mu_M}{\moy{\muMi}{f^2}\bcov{\muMi}{f_i^2}{\sum_{k\in I_i} F'(x_k)}},
  $$
  where $f_i^2 :=f^2/\moy{\muMi}{f^2}$. Thus, by the Cauchy-Schwarz
  inequality, the first term of the RHS of \eqref{in:ceg-step-1} is
  bounded above by
  $$
  2\,\ell\,\sum_{i=1}^\ell
  \bmoy{{\mu_M}}
  {\moy{\muMi}{f^2}\bcov{\muMi}{f_i^2}{\sum_{k\in I_i} F'(x_k)}^2},
  $$
  where we used the Jensen inequality with respect to the density
  $\moy{\muMi}{f^2}$.  Now, by the entropy inequality and by
  \cite[Lemma 6.1]{landim-panizo-yau}
  $$
  \moy{\muMi}{f^2} 
  \bcov{\muMi}{f_i^2}{\sum_{k\in I_i} F'(x_k)}^2
  \leq C\,\ABS{I_i}\,\ent{\muMi}{f^2}.
  $$
  At this stage, the argument used for the Poincar\'e inequality can
  be rewritten exactly in the same way, by replacing the variance by
  the entropy and $\cP_K$ by $\cL_K$. It gives finally that the first
  term of the RHS of \eqref{in:ceg-step-1} is bounded above by
  $$
  nC_K\cL_K\,\moy{\mu_{M}}{\ABS{\GR f}^2}.
  $$
  The latter can be rewritten by virtue of \eqref{eq:mumeqnum} as
  $n\,C'_\veps\,\moy{{\si_M}}{\ABS{\GR f}^2}$.  It remains to bound
  the last term of the RHS of \eqref{in:ceg-step-1}.  Let $\be_0$ as
  in \cite[Lemma 6.5]{landim-panizo-yau} and $\de\in(0,2)$.  By a
  simple rewriting of \cite[Lemma 4.5]{landim-panizo-yau}, one gets
  that if $\ent{\mu_M}{f^2} \leq \de(n+1)\be_0^2$ with $n$ and $K$
  large enough
  $$
  \bcov{{\mu_M}}
  {f^2}{\sum_{i=1}^\ell\ABS{I_i}\moy{\muMi}{F'}}^2
  \leq \de n C\,\ent{\mu_M}{f^2}.
  $$
  In the other hand, if $\ent{\mu_M}{f^2} \geq \de(n+1)\be_0^2$, one gets
  $$
  \bcov{{\mu_M}} {f^2}{\sum_{i=1}^\ell\ABS{I_i}\moy{\muMi}{F'}}^2 \leq \de n
  C_K\,\ent{\mu_M}{f^2} + C_{K,\de}\,n\,\moy{\mu_M}{\ABS{\GR f}^2}.
  $$
  This last estimate is based on a simple rewriting of \cite[Lemma
  4.5]{landim-panizo-yau} together with the following straightforward but
  essential version of \cite[Lemma 4.6]{landim-panizo-yau} :
  $$
  \moy{\nu_{I_i\cup I_j,M}}{(m_i-m_j)^2f^2} \leq C_1(K)\,\moy{\nu_{I_i\cup
      I_j,M}}{f^2} + C_2(K)\,\cL_{2K}\,\bmoy{\nu_{I_i\cup I_j,M}}{\sum_{k\in
      I_i\cup I_j} \ABS{\pd_k f}^2},
  $$
  where $\nu_{I_i\cup I_j,M}$ is the conditional measure on
  $I_i\cup I_j$, $m_i=\ABS{I_i}^{-1}\sum_{k\in I_i}$, and $C_1(K)\to
  0$ when $K\to+\infty$.
  
  Summarising, we get that for any $\de\in (0,2)$ and for $n$ and $K$
  large enough, the last term of the RHS of \eqref{in:ceg-step-1} is
  bounded above as follows
  $$
  \bcov{{\mu_M}} {f^2}{\sum_{i=1}^\ell\ABS{I_i}\moy{\muMi}{F'}}^2 
  \leq %
  n\de C_K\,\ent{{\mu_M}}{f^2}
  +nC_{K,\de}\cL_{2K}\,\moy{\mu_M}{\ABS{{\GR}f}^2},
  $$
  which can be rewritten by virtue of \eqref{eq:mumeqnum} as $\veps
  n C'\,\ent{{\si_M}}{f^2} + nC''_\veps\,\moy{\si_M}{\ABS{\GR f}^2}$.
  This achieves the proof of \eqref{in:ceg} and Lemma \ref{le:ceg}.
\end{proof}

\section*{Acknowledgements}

The author would like to warmly acknowledge Prof. Michel Ledoux for
helpful discussions and encouraging comments, and Doct. Ivan Gentil
for some discussions at the beginning of this work.

\providecommand{\etalchar}[1]{$^{#1}$}
\providecommand{\bysame}{\leavevmode\hbox to3em{\hrulefill}\thinspace}
\providecommand{\MR}{\relax\ifhmode\unskip\space\fi MR }
\providecommand{\MRhref}[2]{%
  \href{http://www.ams.org/mathscinet-getitem?mr=#1}{#2}
}
\providecommand{\href}[2]{#2}

\begin{center}
  \hrule
\end{center}

\noindent
\textbf{E-mail}: \verb+mailto:chafai@math.ups-tlse.fr+\\
\textbf{URL} : \verb+http://www.lsp.ups-tlse.fr/Chafai/+\\
\textbf{Postal Address}: Laboratoire de Statistique et Probabilit\'es, U.M.R.
C.N.R.S. C5583, D\'epartement de Math\'ematiques, Universit\'e Paul Sabatier,\\
118 route de Narbonne, F-31062 Toulouse \textsc{Cedex} 4, France.\\
\textbf{Keywords}: \mykeywords.\\
\textbf{Subj. Class. MSC-2000} : \mysubjclass.


\begin{thebibliography}{ABC{\etalchar{+}}00}

\bibitem[ABC{\etalchar{+}}00]{our-ls-2000}
C.~An\'e, S.~Blach\`ere, D.~Chafa\"{\i}, P.~Foug\`eres, I.~Gentil, F.~Malrieu, C.~Roberto,
  and G.~Scheffer, \emph{Sur les in\'egalit\'es de {S}obolev logarithmiques},
  Panoramas et Synth\`eses, vol.~10, Soci\'et\'e {M}ath\'ematique de {F}rance, Paris,
  2000.

\bibitem[BH99]{hel-bod}
T.~Bodineau and B.~Helffer, \emph{The log-{S}obolev inequality for unbounded
  spin systems}, J. Funct. Anal. \textbf{166} (1999), no.~1, 168--178.

\bibitem[Cap01]{caputo}
P.~Caputo, \emph{A remark on spectral gap and logarithmic sobolev inequalities
  for conservative spin systems}, Preprint (2001) archived as \textbf{mp\_arc 01-71}.

\bibitem[Cap02]{caputo2}
P.~Caputo, \emph{Uniform Poincar\'e inequalities for unbounded conservative spin
systems: The non-interacting case}, Preprint (2002) archived as \textbf{mp\_arc 02-46}.

\bibitem[He98]{helffer-98}
B.~Helffer, \emph{Remarks on decay of correlations and {W}itten
{L}aplacians, {B}rascamp-{L}ieb inequalities and semiclassical limit},
J. Funct. Anal., \textbf{155} (1998), no.~2, 571--586.

\bibitem[He99]{helffer-99}
B.~Helffer, \emph{Remarks on decay of correlations and {W}itten
{L}aplacians. {I}{I}. {A}nalysis of the dependence on the interaction},
Rev. Math. Phys., \textbf{11} (1999), no.~3, 321--336.

\bibitem[He99-2]{helffer-99-2}
B.~Helffer, \emph{Remarks on decay of correlations and {W}itten
{L}aplacians. {I}{I}{I}. {A}pplication to logarithmic {S}obolev inequalities},
Ann. Inst. H. Poincaré Probab. Statist., \textbf{35} (1999), no.~4, 483--508.

\bibitem[KL99]{kipnis-landim}
C.~Kipnis and C.~Landim, \emph{Scaling limit of interacting particle systems},
  Grundlheren der {M}athematischen {W}issenschaften, vol. 320, Springer,
  Berlin, New-York, 1999.

\bibitem[Led01]{ledoux-lsrev}
M.~Ledoux, \emph{Logarithmic {S}obolev inequalities for unbounded spin systems
  revisited}, S\'eminaire de Probabilit\'es, XXXV, Springer, Berlin, 2001,
  pp.~167--194.

\bibitem[LPY00]{landim-panizo-yau}
C.~Landim, G.~Panizo, and H.~T. Yau, \emph{Spectral gap and logarithmic
  {S}obolev inequality for unbounded conservative spin systems}. 
  Ann. Inst. H. Poincar{\'e} PR \textbf{38} (2002), no.~5, pp.~739--777.

\bibitem[LY93]{lu-yau}
S.~L. Lu and H.-T. Yau, \emph{Spectral gap and logarithmic {S}obolev inequality
  for {K}awasaki and {G}lauber dynamics}, Comm. Math. Phys. \textbf{156}
  (1993), no.~2, 399--433.

\bibitem[Mal01]{malrieu-projection}
F.~Malrieu, \emph{Convergence to equilibrium for granular media equations and
  their {E}uler schemes}, preprint 2001, 
  \verb+http://www.lsp.usp-tlse.fr/Fp/Malrieu/+,
  to appear in Annals of Applied Probability.

\bibitem[Mar99]{martinelli-99}
F.~Martinelli, \emph{Lectures on {G}lauber dynamics for discrete spin models},
  Lectures on probability theory and statistics. {\'E}cole d'{\'e}t{\'e} de
  probabilit{\'e}s de St-Flour 1997, Lecture Notes in Math., vol. 1717,
  Springer, Berlin, 1999, pp.~93--191.

\bibitem[SC97]{saloffcoste}
L.~Saloff-Coste, \emph{Lectures on finite {M}arkov chains}, Lectures on
  probability theory and statistics. {\'E}cole d'\'et\'e de probabilit\'es de St-Flour
  1996, Lecture Notes in Math., vol. 1665, Springer, Berlin, 1997,
  pp.~301--413.

\bibitem[SZ92]{stroockzegarlinski-92-b}
D.~Stroock and B.~Zegarli{\'n}ski, \emph{The logarithmic {S}obolev inequality
  for continuous spin systems on a lattice}, J. Funct. Anal. \textbf{104}
  (1992), no.~2, 299--326.

\end{thebibliography}
\end{document}